\documentclass[12pt,a4wide]{article}

\usepackage{a4wide}

\usepackage[intlimits]{amsmath}
\usepackage{amssymb,amsthm}

\usepackage[cp1251]{inputenc}
\usepackage[english,russian]{babel}
\begin{document}
\vskip 20mm {\large

\hskip 5mm\textbf{On weak convergence for stochastic evolutionary
systems}

\hskip 5cm\textbf{in average principle.}

\vskip 1cm\hskip 55mm\textbf{I.V.Samoilenko\\
Institute of Mathematics,\\ Ukrainian National Academy of
Sciences,\\ 3 Tereshchenkivs'ka, Kyiv, 01601, Ukraine\\
isamoil@imath.kiev.ua}

\vskip 1cm
 \textbf{Abstract:}\textit{Weak convergence of the
stochastic evolutionary system to the average evolutionary system is
proved. The method proposed by R.Liptser in [4] for semimartingales
is used.
 But we apply a solution of singular perturbation problem instead of ergodic theorem.}

\vskip 1cm

Many authors studied the problems of weak convergence for different
classes of stochastic processes (see [1-7]). The methods used are
also different and depend on the process studied: ergodic theorems,
relative compactness, martingales, etc.

We propose to study a stochastic evolutionary system by a
combination of two methods. The method proposed by R.Liptser in [4]
for semimartingales is combined with a solution of singular
perturbation problem instead of ergodic theorem. We use

\textbf{Theorem [1].} \textit{Let the following conditions hold for
a family of Markov processes $\xi^{\varepsilon}(t), t\geq0,
\varepsilon>0$:}

\textit{\textbf{C1}: There exists a family of test functions
$\varphi^{\varepsilon}(u, x)$ in $C^2_0(\mathbb{R}^d\times E)$, such
that $$\lim\limits_{\varepsilon\to 0}\varphi^{\varepsilon}(u, x) =
\varphi(u),$$ uniformly on $u, x.$}

\textit{\textbf{C2}: The following convergence holds
$$\lim\limits_{\varepsilon\to
0}\mathbb{L}^{\varepsilon}\varphi^{\varepsilon}(u, x) =
\mathbb{L}\varphi(u),$$ uniformly on $u, x$. The family of functions
$\mathbb{L}^{\varepsilon}\varphi^{\varepsilon}, \varepsilon>0$ is
uniformly bounded, and $\mathbb{L}\varphi(u)$ and
$\mathbb{L}^{\varepsilon}\varphi^{\varepsilon}$ belong to
$C(\mathbb{R}^d\times E)$.}

\textit{\textbf{C3}: The quadratic characteristics of the
martingales that characterize a coupled Markov process
$\xi^{\varepsilon}(t), x^{\varepsilon}(t), t\geq0, \varepsilon>0$
have the representation $<\mu^{\varepsilon}>_t = \int^t_0
\zeta^{\varepsilon}(s)ds,$ where the random functions
$\zeta^{\varepsilon}, \varepsilon> 0,$ satisfy the condition
$$\sup\limits_{0\leq s \leq T} \mathbf{E}|\zeta^{\varepsilon}(s)|\leq
c < +\infty.$$}

\textit{\textbf{C4}: The convergence of the initial values holds and
$$\sup\limits_{\varepsilon>0}\mathbf{E}|\zeta^{\varepsilon}(0)|\leq C
< +\infty.$$}

\textit{Then the weak convergence $$\xi^{\varepsilon}(t)\Rightarrow
\xi(t), \varepsilon\to 0,$$ takes place.}

The evolutionary system we are going to study is constructed in such
a way that the conditions \textbf{C1, C2} are satisfied as a result
of solution of singularly perturbation problem for the system. To
verify the conditions \textbf{C3, C4} we use the method proposed in
[4] extending Bogolubov averaging principle for discontinuous
semimartingales.

Really, we study a stochastic evolutionary system with the switching
ergodic Markov process
$$du^{\varepsilon}_t/dt=b(u^{\varepsilon}_t;{\ae}(t/\varepsilon)), \eqno(1)$$
$$u^{\varepsilon}_0=u.$$

The switching Markov process ${\ae}(t), t\geq0$ on the standard
(Polish) phase space $(E,\mathcal{E})$ is defined by the generator
$$Q\varphi(x)=\int_E Q(x,dy)[\varphi(y)-\varphi(x)],$$ where
$$Q(x,dy)=q(x)P(x,dy).$$

The stationary distribution $\pi(dx)$ of the egrodic process defines
the projector $$\Pi\varphi(x)=\widehat{\varphi}\mathbf{1}(x),
\widehat{\varphi}:=\int_E\pi(dx)\varphi(x), \mathbf{1}(x)\equiv 1.$$

We also define a potential $R_0$ that is, for ergodic Markov process
with generator $Q$ and semigroup $P_t, t\geq0$, a bounded operator
$$R_0=\int_0^{\infty}(P_t-\Pi)dt,$$ $$QR_0=R_0Q=\Pi-I.$$

The velocity function $b(u;x), u\in \mathbb{R}^d, x\in E$ provides a
global solution of the deterministic equations
$$du_x(t)/dt=b(u_x(t);x), u_x(0)=u, x\in E.$$

Markov process $u^{\varepsilon}_t, {\ae}^{\varepsilon}(t):={\ae}(t/
\varepsilon), t\geq0$ can be characterized by the generator [1,
Proposition 3.3]
$$\mathbb{L}^{\varepsilon}\varphi(u;x)=[\varepsilon^{-1}Q+\mathbb{B}(x)]\varphi(u;x), \eqno(2)$$
where the generator $\mathbb{B}(x)\varphi(u)=b(u;x)\varphi'(u).$

It follows from Proposition 5.6 [1] that the solution of the
singular perturbation problem for the generator (2) is given by the
relation
$$\mathbb{L}^{\varepsilon}\varphi^{\varepsilon}(u;x)=\widehat{\mathbb{B}}\varphi(u)
+\varepsilon\theta^{\varepsilon}(x)\varphi(u) \eqno(3)$$ on the
functions
$\varphi^{\varepsilon}(u;x)=\varphi(u)+\varepsilon\varphi_1(u;x)$
with $\varphi(u) \in C_0^2(\mathbb{R}^d).$

Thus, we see from (3) that the solution of singularly perturbation
problem for $\mathbb{L}^{\varepsilon}, \varphi^{\varepsilon}(u;x)$
satisfies the conditions \textbf{C1, C2}.

The operator
$$\widehat{\mathbb{B}}\varphi(u)=\widehat{b}(u)\varphi'(u), \widehat{b}(u)=\int_E\pi(dx)b(u;x).$$

The term
$\theta^{\varepsilon}(x)=\mathbb{B}(x)R_0\widetilde{\mathbb{B}}(x),
\widetilde{\mathbb{B}}(x):= \mathbb{B}(x)-\widehat{\mathbb{B}}$
satisfy the condition
$$|\theta^{\varepsilon}\varphi(u)|\leq C<\infty.$$

\textit{Remark.} The definition of the operator
$\widetilde{\mathbb{B}}(x)$ depends on the definition of the
potential $R_0$. Namely, if $QR_0=\Pi-I,$ as we defined, then
$\widetilde{\mathbb{B}}(x):= \mathbb{B}(x)-\widehat{\mathbb{B}}.$ If
we'd define $QR_0=I-\Pi,$ then we'd have
$\widetilde{\mathbb{B}}(x):=\widehat{\mathbb{B}}-\mathbb{B}(x).$

Let us now regard an average evolutionary system
$$d\widehat{u}_t/dt=\widehat{b}(\widehat{u}_t), \widehat{u}_0=u.\eqno(4)$$

Our aim is to prove a weak convergence of the stochastic
evolutionary system (1) to the average evolutionary system (4). To
do this, we should show that the conditions \textbf{C1-C4} of the
Theorem are satisfied. As we showed, the conditions \textbf{C1, C2}
are true due to the solution of singularly perturbation problem for
the evolutionary system. Thus, we are going to verify \textbf{C3,
C4}.

Following [4] we demand from the function $b(u;x)$ to satisfy the
the following conditions:

$\mathbf{I}$ (linear growth) $$|b(u;x)|\leq L(1+|u|),$$

$\mathbf{II}$ (Lipschitz continuity) $$|b(u;x)-b(u';x)|\leq
C|u-u'|.$$

To prove the main result we need the following lemma.

{\bf Lemma.} Under the assumption $\mathbf{I}$ there exists a
constant $k>0$, independent of $\varepsilon$ and dependent on $T$:
$$\mathbf{E}\sup\limits_{t\leq T}|u^{\varepsilon}_t|^2\leq k_T.$$

{\bf Corollary.} (follows from Chebyshev's inequality) Under the
assumption $\mathbf{I}$ compact containment condition (CCC) holds:
$$\lim\limits_{c\to \infty}\sup\limits_{\varepsilon\leq\varepsilon_0}
\mathbf{P}\{\sup\limits_{t\leq T}|u^{\varepsilon}_t|>c\}=0.$$

{\it Proof of Lemma:} (follows [4]) Let us rewrite the equation (1)
in the form:
$$u^{\varepsilon}_t=u+\int_0^tb(u^{\varepsilon}_s;{\ae}(s/\varepsilon))ds=:u+A_t^{\varepsilon}.$$

If we put $u_{t}^*=\sup\limits_{s\leq t}|u_t|,$ then
$$((u^{\varepsilon}_t)^*)^2\leq 2[u^2+((A^{\varepsilon}_t)^*)^2].\eqno(6)$$

The condition $\mathbf{I}$ implies that
$$(A^{\varepsilon}_t)^*\leq L\int_0^t(1+(u^{\varepsilon}_s)^*)ds.$$

The last inequality with (6) and Cochy-Bounikovsky inequality
([$\int_0^t\varphi(s)ds]^2\leq t\int_0^t\varphi^2(s)ds$) implies,
with some constants $k_1$ and $k_2$, independent of $\varepsilon$,
that
$$\mathbf{E}((u^{\varepsilon}_t)^*)^2\leq
k_1+k_2\int_0^t\mathbf{E}((u^{\varepsilon}_s)^*)^2ds,$$ for $t\leq
T.$

By Gronwall-Bellman inequality we obtain
$$\mathbf{E}((u^{\varepsilon}_t)^*)^2\leq k_1\exp(k_2 T).$$

Lemma is proved.

{\it Remark.} Another way to prove CCC is proposed in [1, Theorem
8.10] and used by other authors [2,3]. They use the function
$\varphi(u)=\sqrt{1+u^2}$ and prove corollary for
$\varphi(u^{\varepsilon}_t)$ by applying the martingale
characterization of the Markov process.

This may be easily done, due to specific properties of $\varphi(u).$







The main result of our work is the following.

\textbf{Proposition.}  Under the conditions \textbf{I, II}
$$\mathbf{P} - \lim\limits_{\varepsilon\to 0}\sup\limits_{t\leq
T}|u^{\varepsilon}_t-\widehat{u}_t|=0.\eqno(5)$$

{\it Proof of Proposition:} As we mentioned, to prove the
Proposition we should verify the conditions \textbf{C3, C4} of the
Theorem.

\textit{Remark.} Condition \textbf{C3} of the Theorem means that the
quadratic characteristics of the martingale, corresponding to a
coupled Markov process, is relatively compact. The same result
follows from the CCC (see Corollary) by [6].

Thus, the condition \textbf{C3} follows from the Corollary.





As soon as $u^{\varepsilon}_0=\widehat{u}_0=u,$ we see from the
Lemma that the condition \textbf{C4} follows from \textbf{I, II}.

Thus, all the conditions of the Theorem 6.3 [1] are satisfied, so
(5) is true.

Proposition is proved.

{\it Remark.} To compare the results obtained in this work and in
the work [4], we should note that R.Liptser uses the following
scheme: he first applies the conditions $\textbf{I}, \textbf{II}$ to
obtain the result of Lemma 1, and CCC follows from this Lemma. We
also apply the conditions $\textbf{I}, \textbf{II}$ to obtain Lemma
1 and CCC.

Then, by [4], the necessary weak convergence is obtained by the use
of ergodic theorem. But we, on the contrary, use singularly
perturbation problem, and so the weak convergence follows from the
conditions $\mathbf{C1, C2}$.

Similar method is used by A.V.Skorokhod in [7,8]. But his condition
of weak convergence looks like existence of a generator, equivalent
to our $\widehat{\mathbb{B}}$, that corresponds to a functional of
the process. This generator should approximate the functional's
surplus.

Another method is used in [2]. Here a well-known conditions of
compact- ness are applied to the family of the processes studied.
Then the weak convergence follows from compactness and convergence
at a "good" class of test-functions.

In [3,5] similar conditions of compactness are used for a family of
processes. Then a submartingale, corresponding to the family of
processes is constructed. The weak convergence follows then from the
compactness and uniqueness of the solution of martingale problem.












\newpage
\noindent\textbf{\large References}

\begin{enumerate}

 \item  \textbf{Koroliuk, V.S., Limnios, N.} Stochastic Systems in Merging Phase Space.
\textit{ World Scientific Publishers (2005), 330p.}

 \item  \textbf{Ethier, S.N., Kurtz T.G.} Markov Processes: Characterization
and convergence. \textit{ J. Wiley $\&$ Sons, New York, (1986),
529p.}

 \item  \textbf{Свириденко, М.Н.}  Мартингальный подход к получению предельных теорем для марковских процессов.
     \textit{ВИНИТИ №37, (1986), 30с.}

  \item \textbf{Liptser, R.Sh.}  The Bogolubov averaging principle for semimartingales.
     \textit{Proceedings of the Steclov Institute of Mathematics, No.4, (1994), 12p.}

\item \textbf{Stroock, D.W., Varadhan, S.R.S.} Multidimensional
Diffusion Processes. \textit{Springer-Verlag, Berlin, (1979), 339p.}

\item \textbf{Jacod J., Shiryaev A.N.} Limit Theorems for Stochastic
Processes. \textit{Springer-Verlang (1987), Berlin, 325p.}

\item \textbf{Скороход, А.В.} Асимптотические методы теории
стохастических дифференциальных уравнений. \textit{Наукова думка,
Киев, (1987), 328с.}

\item \textbf{Skorokhod, A.V., Hoppensteadt, F.C., Salehi, H.} Random
Perturbation Methods with Applications in Science and Engineering.
\textit{Springer (2002), 488p.}

\end{enumerate}

\end{document}